\documentclass[11pt]{article}
\usepackage{amsmath}
\usepackage{amssymb}
\makeatletter
\newif\ifmsbmloaded@
\def\loadmsbm{\msbmloaded@true
  \font\tenmsb=msbm10 scaled 1\@ptsize00
  \font\sevenmsb=msbm7 scaled 1\@ptsize00
  \font\fivemsb=msbm5 scaled 1\@ptsize00
  \alloc@8\fam\chardef\sixt@@n\msbfam
  \textfont\msbfam=\tenmsb
  \scriptfont\msbfam=\sevenmsb
  \scriptscriptfont\msbfam=\fivemsb
  }
\@addtoreset{equation}{section}

\makeatother \loadmsbm
\def\R{\Bbb R}

\def\N{\Bbb N}
\def\Z{\Bbb Z}

\def\ep{\epsilon}
\def\pa{\partial}

\def\th{\theta}

\def\na{\nabla}

\def\al{\alpha}

\def\endproof{\hphantom{MM}\hfill\llap{$\square$}\goodbreak}

\newcommand{\beq}{\begin{equation}}
\newcommand{\eeq}{\end{equation}}
\newcommand{\ben}{\begin{eqnarray}}
\newcommand{\een}{\end{eqnarray}}
\newcommand{\beno}{\begin{eqnarray*}}
\newcommand{\eeno}{\end{eqnarray*}}

\newtheorem{Theorem}{Theorem}[section]
\newtheorem{Lemma}[Theorem]{Lemma}
\newtheorem{Def}{Definition}[section]
\newtheorem{Remark}{Remark}[section]
\newtheorem{Proposition}[Theorem]{Proposition}

\begin{document}

\title{A new Bernstein's Inequality and the 2D Dissipative Quasi-Geostrophic Equation}
\date{}
\author{Qionglei Chen $^1$,
Changxing Miao $^2$,
Zhifei Zhang $^3$\\
{\small  $^{1, 2}$Institute of Applied Physics and Computational Mathematics,}\\
    {\small P.O. Box 8009, Beijing 100088, P. R. China.}\\
     {\small (chen\_qionglei@iapcm.ac.cn  and   miao\_changxing@iapcm.ac.cn)}\\
{\small $^{3}$ School of Mathematical Science, Peking University,}\\
{\small Beijing 100871, P. R.  China.}\\
{\small(zfzhang@math.pku.edu.cn)  }}
\maketitle

\begin{abstract}
We show a new Bernstein's inequality which generalizes the results
of Cannone-Planchon, Danchin and Lemari\'{e}-Rieusset. As an
application of this inequality, we prove the global well-posedness
of the 2D quasi-geostrophic equation with the critical and
super-critical dissipation for the small initial data in the
critical Besov space, and local well-posedness for the large
initial data.

\end{abstract}
\vspace*{3mm}


\section{Introduction}
\quad\,
We are concerned with the 2D dissipative
quasi-geostrophic equation
\begin{equation}\label{1.1}
(QG)_\alpha\quad\left\{
\begin{aligned}
&\partial_t\theta+u\cdot \nabla \theta+\kappa(-\Delta)^\alpha
\theta=0,\quad x\in \R^2,\, t>0,\\
&\theta(0,x)=\theta_0(x).
\end{aligned}
\right.
\end{equation}
Here $\alpha\in [0,\frac 1 2], \kappa>0$ is the dissipative coefficient,
$\theta(t,x)$ is a real-valued function of $t$ and $x$. The
function $\theta$ represents the potential temperature, the fluid
velocity $u$ is determined from $\theta$ by a stream function
$\psi$
\begin{equation}\label{1.2}
(u_1,u_2)=\bigg(-\frac {\partial \psi} {\partial x_2},\frac {\partial
\psi}{\partial x_1}\bigg), \quad (-\Delta)^\frac 1 2\psi=-\theta.
\end{equation}
A fractional power of the Laplacian  $(-\Delta)^\beta$ is defined
by
$$\widehat{(-\Delta)^\beta f}(\xi)=|\xi|^{2\beta}\hat f(\xi),$$
where $\hat f$ denotes the Fourier transform of $f$. We rewrite (\ref{1.2}) as
\begin{equation}
u=(\partial_{x_2}(-\Delta)^{-\frac 1 2}\theta,
-\partial_{x_1}(-\Delta)^{-\frac 1 2}\theta)={\cal
R}^\bot\theta=(-{\cal R}_2\theta,{\cal R}_1\theta),\nonumber
\end{equation}
where ${\cal R}_k, k=1,2$, is the Riesz transform defined by
$$\widehat{{\cal R}_kf}(\xi)=-i\xi_k/|\xi|\hat f(\xi).$$

$(QG)_\alpha$ is an important model in geophysical fluid dynamics,
they are special cases of the general quasi-geostrophic
approximations for atmospheric and oceanic fluid flow with small
Rossby and Ekman numbers.  There exists deep analogy between the
equation (\ref{1.1}) with $\al=\frac{1}{2}$ and the the 3D
Navier-Stokes equations. For more details about its  background in
geophysics, see \cite{Con1,Ped}. The case $\alpha>\frac 1 2$ is
called the subcritical case, the case $\alpha=\frac 1 2$ is
critical , and the case $0\le \alpha<\frac 1 2$ is  supercritical.
In the subcritical case, Constantin and Wu\cite{Con2} proved the
existence of global in time smooth solutions. In the critical
case, Constantin, Cordoba, and Wu\cite{Con3} proved the existence
and uniqueness of global smooth solution on the spatial periodic
domain under the assumption of small $L^\infty$ norm. Recently,
Chae and Lee\cite{Cha2} studied the super-critical case and proved
the global well-posedness  for small data in the Besov spaces
$\dot B^{2-2\alpha}_{2,1}$. Very recently,
Cordoba-Cordoba\cite{Cor3}, Ning\cite{Ju1,Ju2} studied the
existence and uniqueness in the Sobolev spaces $H^s, s\ge
2-2\alpha, \alpha\in [0,\frac 1 2]$. Wu \cite{Wu1,Wu2} studied the
well-posedness in general Besov space $B^s_{p,q}, s>2(1-\alpha),
p=2^N$.  Many other relevant results can also be found in
\cite{Cha1,Cor1,Cor2}.

One purpose of this paper is to study the well-posedness of the 2D dissipative
quasi-geostrophic equation in the critical Besov space
$B^{\frac 2 p+1-2\alpha}_{p,q}, p\ge 2, q\in [1,\infty)$.
If we use the standard energy method as in \cite{Cha2, Wu3},
we need to establish the lower bound for the term generated from the dissipative part
\begin{align}\label{1.3}
\int_{\R^2}\Lambda^{2\alpha}\Delta_j\theta
|\Delta_j\theta|^{p-2}\Delta_j\theta dx\ge 2^{2\al j}\|\Delta_j\theta \|_p^p, \quad p\ge 2,
\end{align}
where $\Delta_j$ is the frequency localization operator at
$|\xi|\approx 2^j$ (see Section 2).
For $p=2$, this is a direct consequence of Plancherel formula. In the case $\alpha=1$,
it is proved by Cannone and Planchon\cite{CP}.
To generalize (\ref{1.3}) to general index $\alpha, p$, it is sufficient to show the following Bernstein's inequality
\begin{equation}\label{1.4}
c_p2^{\frac{2\alpha j}{p}}\|\Delta_jf\|_p\le
\|\Lambda^\alpha(|\Delta_jf|^{\frac{p}{2}})\|_2^{\frac{2}{p}}
\le C_p2^{\frac{2\alpha j}{p}}\|\Delta_jf\|_p, \quad p>2,
\end{equation}
which together with an improved positivity Lemma 3.3 in \cite{Ju2}
(see also Section 3, Lemma \ref{Lem3.3}) will imply (\ref{1.3}).
We should point out that (\ref{1.4}) is proved by
Lemari\'{e}-Rieusset\cite{LR} in the case $\al=1$, and by Danchin
\cite{Dan1} when $p$ is any even integer. On the other hand, Wu
\cite{Wu3} gives a formal proof for general index. The first
purpose of this paper is to present a rigorous proof of Theorem
3.4 in \cite{Wu3} which plays a key role in Wu's paper.

\begin{Theorem}[Bernstein's inequality]\label{Thm1}
Let $p\in [2,\infty)$ and $\alpha \in [0,1]$. Then there exist
 two positive constants $c_p$ and $C_p$ such that
for any $f\in\cal{S}'$ and $j\in\Z$, we have
\begin{equation}\label{1.6}
c_p2^{\frac{2\alpha j}{p}}\|\Delta_jf\|_p\le
\|\Lambda^\alpha(|\Delta_jf|^{\frac{p}{2}})\|_2^{\frac{2}{p}}
\le C_p2^{\frac{2\alpha j}{p}}\|\Delta_jf\|_p.
\end{equation}
\end{Theorem}

The second purpose is to study the well-posedness of the 2D
dissipative quasi-geostrophic equation in the critical Besov space
$B^{\frac 2 p+1-2\alpha}_{p,q}$ by  using Theorem \ref{Thm1} and
Fourier localization technique.

\begin{Theorem}\label{Thm2}
Assume that $(\alpha,p,q)\in(0,\frac 1 2]\times[2,\infty)\times [1,\infty)$.  If
$\theta_0$ belongs to ${B^\sigma_{p,q}}$  with
$\sigma=\frac 2 p+1-2\alpha,$ then
there exists a positive real number $T$ such that a unique solution
to the 2D dissipative quasi-geostrophic equation
$\theta(t,x)$ exists on $[0,T)\times \R^2$
satisfying
$$
\theta(t,x) \in C([0,T); B^\sigma_{p,q})\cap \widetilde{L}^1(0,T;\dot{B}^{\frac{2}{p}+1}_{p,q}),
$$
with the time $T$ bounded from below by
$$\sup\big\{T'>0: \,\big\|(1-e^{-\kappa c_p2^{2\al j}T'})^\frac{1}{2}2^{j\sigma}
\|\Delta_j\theta_0\|_p\big\|_{\ell^q(\Z)}
\le c\kappa\big\}.$$
Furthermore, if $\|\theta_0\|_{\dot B^\sigma_{p,q}}\le\epsilon\kappa$ for some positive number $\epsilon$,
then we can choose $T=+\infty$.
\end{Theorem}

\begin{Remark}
It is pointed out that the homogeneous Besov  space $\dot B^{\sigma}_{p,q}$ is important as it
gives the important scaling invariant function space. In fact, if
$\theta(t,x)$ and $u(t,x)$ are solutions of (\ref{1.1}), then
$\theta_\lambda(t,x)=\lambda^{2\alpha-1}\theta(\lambda^{2\alpha
}t, \lambda x)$ and $u_\lambda(t,x)=\lambda^{2\alpha-1}u(\lambda^{2\alpha }t, \lambda
x)$ are also solutions of (\ref{1.1}). The $\dot B^{\sigma}_{p,q}$ norm of $\theta(t,x)$ is invariant under this scaling.
Moreover, for the global existence result, the smallness assumption is imposed only on the homogenous
norm of the initial data.
\end{Remark}

\begin{Remark} The result of Theorem \ref{Thm2} for the case $(p,q)=(2,1)$ corresponds to the result
of Chae and Lee\cite{Cha2} in the critical Besov space $B^{2-2\alpha}_{2,1}$. In the case $(p,q)=(2,2)$ , it
corresponds to the result of Ning\cite{Ju1} in the Sobolev space $H^{2-2\alpha}$.
On the other hand, thanks to the embedding relationship:
$$B^s_{2,1}\subsetneq H^s\subsetneq B^s_{2,q},
\quad \textrm{for} \,\, q>2,$$
our result  improves the results of \cite{Cha2} and \cite{Ju1}.
\end{Remark}

\begin{Remark}
Wu\cite{Wu1,Wu2} proved the well-posedness of (\ref{1.1}) for the initial data
in the sub-critical Besov space $B^s_{p,q}$ with $s>2-2\alpha, p=2^N$.
We obtain the well-posedness in the critical Besov space $B^{\frac 2 p+1-2\alpha}_{p,q}$,
and get rid of the restriction on $p=2^N$.
\end{Remark}

\begin{Remark} Very recently, Miura\cite{Miu} proved the local well-posedness of (\ref{1.1})
for the large initial data in the critical Sobolev space
$H^{2-2\alpha}$. His result is a particular case of Theorem
\ref{Thm2}, and our proof is  simpler(see section 4.2).
\end{Remark}

\noindent{\bf Notation:}\,\, Throughout the paper, $C$ denotes various ``harmless'' large
finite constants, and $c$  denotes various ``harmless" small
constants. We shall sometimes use $X\lesssim Y$ to
denote the estimate $X\le CY$ for some $C$.
$\{c_j\}_{j\in\Z}$ denotes any positive series with $\ell^q(\Z)$ norm less than or equals to 1.
we shall sometimes use the $\|\cdot\|_p$ to
denote $L^p(\R^d)$ norm of a function.

\section{Littlewood-Paley decomposition}

\quad\, Let us recall the Littlewood-Paley decomposition. Let
${\cal S}(\R^d)$ be the Schwartz class of rapidly decreasing
functions. Given $f\in {\cal S}(\R^d)$, its Fourier transform
${\cal F}f=\hat f$ is defined by
$$
\hat f(\xi)=(2\pi)^{-\frac{d}2}\int_{\R^d}e^{-ix\cdot \xi}f(x)dx.
$$
Choose two nonnegative radial functions $\chi$, $\varphi \in {\cal
S}(\R^d)$, supported respectively in ${\cal B}=\{\xi\in\R^d,\,
|\xi|\le\frac{4}{3}\}$ and ${\cal C}=\{\xi\in\R^d,\,
\frac{3}{4}\le|\xi|\le\frac{8}{3}\}$ such that
\beno
\chi(\xi)+\sum_{j\ge0}\varphi(2^{-j}\xi)=1,\quad\xi\in\R^d,\\
\sum_{j\in\Z}\varphi(2^{-j}\xi)=1,\quad\xi\in\R^d\backslash \{0\}.
\eeno Setting $\varphi_j(\xi)=\varphi(2^{-j}\xi)$. Let $h={\cal
F}^{-1}\varphi$ and $\tilde{h}={\cal F}^{-1}\chi$, we  define the
frequency localization operator as follows \beno
&&\Delta_jf=\varphi(2^{-j}D)f=2^{jd}\int_{\R^d}h(2^jy)f(x-y)dy, \\
&&S_jf=\sum_{k\le j-1}\Delta_kf=\chi(2^{-j}D)f=2^{jd}\int_{\R^d}\tilde{h}(2^jy)f(x-y)dy.
\eeno
Informally, $\Delta_j=S_{
j}-S_{j-1}$  is a frequency projection to the annulus
$\{|\xi|\approx 2^j\}$, while $S_j$ is a frequency projection to the
ball $\{|\xi|\lesssim 2^j\}$. One easily verifies that
with our choice of $\varphi$
\begin{eqnarray}\label{2.1}
\Delta_j\Delta_kf\equiv0\quad i\!f\quad|j-k|\ge 2\quad and
\quad \Delta_j(S_{k-1}f\Delta_k
f)\equiv0\quad i\!f\quad|j-k|\ge 5.
\end{eqnarray}

Now we give the definitions of the Besov spaces.
\begin{Def}\label{Def2.1}Let
$s\in \R, 1\le p,q\le\infty$, the homogenous Besov space $\dot
{B}^s_{p,q}$ is defined by
$$\dot {B}^s_{p,q}=\{f\in {\cal Z}'(\R^d); \|f\|_{\dot
{B}^s_{p,q}}<\infty\}.$$ \end{Def}Here
$$\|f\|_{\dot{B}^s_{p,q}}=\left\{\begin{array}{l}
\displaystyle\bigg(\sum_{j\in\Z}2^{jsq}\|\Delta_j f\|_p^q\bigg)^{\frac 1
q},\quad \hbox{for}\quad q<\infty,\\
\displaystyle\sup_{j\in \Z}\|\Delta_jf\|_p, \quad \hbox{ for} \quad q=\infty,
\end{array}\right.
$$
and ${\cal Z}'(\R^d)$ denotes the dual space of ${\cal Z}(\R^d)=\{
f\in {\cal S}(\R^d); \partial^\gamma\hat f(0)=0; \forall \gamma\in \N^d
\,\hbox {multi-index}\}$ and can be identified by the quotient
space of ${\cal S}'/{\cal P}$ with the polynomials space ${\cal
P}$.

\begin{Def}\label{Def2.1}Let
$s\in \R, 1\le p,q\le\infty$, the inhomogenous Besov space $
{B}^s_{p,q}$ is defined by
$${B}^s_{p,q}=\{f\in {\cal S}'(\R^d); \|f\|_{
{B}^s_{p,q}}<\infty\}.$$ \end{Def}Here
$$\|f\|_{{B}^s_{p,q}}=\left\{\begin{array}{l}
\displaystyle\bigg(\sum_{j\ge 0}2^{jsq}\|\Delta_j f\|_p^q\bigg)^{\frac 1
q}+\|S_0(f)\|_p,\quad \hbox{for}\quad q<\infty,\\
\displaystyle\sup_{j\ge 0}\|\Delta_jf\|_p+\|S_0(f)\|_p, \quad \hbox{ for} \quad q=\infty.
\end{array}\right.
$$
If $s>0$, then ${B}^s_{p,q}=L^p\cap\dot{B}^s_{p,q}$ and $\|f\|_{B^s_{p,q}}\thickapprox\|f\|_{p}+\|f\|_{\dot{B}^s_{p,q}}.$
We refer to \cite{Ber,Tri} for  more details.

Next let's recall Chemin-Lerner's space-time space which will play an
important role in the proof of Theorem \ref{Thm2}.

\begin{Def}\label{2.2}
Let $s\in \R,$ $1\le p, q, r\le\infty$, $I\subset
\R$ is an interval. The homogeneous mixed time-space Besov
space $\widetilde{L}^r(I; \dot B^s_{p,q})$ is the space of the
distribution such that
$$\widetilde{L}^r(I; \dot B^s_{p,q})=
\{f\in {\cal D}(I; {\cal Z'}(\R^{d}));\,\|f\|_{\widetilde{L}^r(I; \dot{B}^{s}_{p,r})}<+\infty\}.$$
\end{Def}
Here
$$
\|f(t)\|_{\widetilde{L}^r(I;\dot B^s_{p,q})}=
\displaystyle\bigg\|2^{sj}
\bigg(\int_{I}\|\Delta_jf(\tau)\|_{p}^rd\tau\bigg)^{\frac{1}{r}}\bigg\|_{\ell^q(\Z)},\quad
$$(usual modification if $r,q=\infty$).
We also need the inhomogeneous mixed time-space Besov space
$\widetilde{L}^r(I; B^s_{p,q})$,  $s>0$, whose norm is defined by
$$
\|f(t)\|_{\widetilde{L}^r(I;
B^s_{p,q})}=\|f(t)\|_{L^r(I;L_x^p)}+\|f(t)\|_{\widetilde{L}^r(I;
\dot B^s_{p,q})}
$$
For the convenience, we sometimes use $\widetilde{L}^r_T(\dot
B^s_{p,q})$ and $\widetilde{L}^r(\dot B^s_{p,q})$ to denote
$\widetilde{L}^r(0,T;\dot B^s_{p,q})$ and
$\widetilde{L}^r(0,\infty;\dot B^s_{p,q})$, respectively.
The direct consequence of Minkowski's inequality is that
$$L^r_t(\dot{B}^{s}_{p,q})\subseteq\widetilde{L}^r_t(\dot{B}^{s}_{p,q})
\quad\hbox{if}\quad r\le q \quad\textrm{and}\quad
\widetilde{L}^r_t(\dot{B}^{s}_{p,q})\subseteq{L}^r_t(\dot{B}^{s}_{p,q})
\quad\hbox{if}\quad r\ge q.$$
We refer to \cite{Che} for more details.\\

Let us state some basic properties about the Besov spaces.
\begin{Proposition}\label{Prop2.1}
$(\rm{i})$ We have the equivalence of norms
$$\|D^k f\|_{\dot B^s_{p,q}}\sim \|f\|_{\dot B^{s+k}_{p,q}}, \quad \textrm{for}\quad  k\in \Z^+.$$
$(\rm{ii})$ Interpolation: for
$s_1, s_2\in\R$ and $\theta\in[0,1]$, one has $$\|f\|_{\dot B^{\theta s_1+(1-\theta)s_2}_{p,q}}\le
\|f\|^\theta_{\dot B^{s_1}_{p,q}}\|f\|^{(1-\theta)}_{\dot B^{s_2}_{p,q}},$$
and the similar interpolation inequality holds for inhomogeneous Besov space.\\
$(\rm{iii})$ Embedding: If $s>\frac{d}{p}$, then
$$B^s_{p,q}\hookrightarrow L^\infty;$$
If  $p_1\le p_2$ and
$s_1-\frac{d}{p_1}>s_2-\frac{d}{p_2}$,  then
$$B^{s_1}_{p_1,q_1}\hookrightarrow B^{s_2}_{p_2,q_2},\quad
B^{s}_{p,\min(p,2)}\hookrightarrow H^{s}_{p}\hookrightarrow B^{s}_{p,\max(p,2)}.$$
Here $H^{s}_{p}$ is the inhomogeneous Sobolev space.
\end{Proposition}
{\it Proof:}\,\,The proof of $(\rm{i})-(\rm{iii})$ is rather standard and one can refer to
\cite{Tri}.\\

Finally we introduce the well-known  Bernstein's inequalities which
will be used repeatedly in this paper.
\begin{Lemma}\label{Lem2.2}
Let $\cal C$ be a ring,  and $\cal B$ a ball, $1\le p\le q\le+\infty$. Assume that $f\in {\cal S'}(\R^d)$,
then for any $|\gamma|\in\Z^+\cup\{0\}$
there exist
constants $C$, independent of $f$, $j$ such that
\begin{align}
&\|\partial^\gamma f\|_q\le C\lambda^{{|\gamma|}+d(\frac{1}{p}-\frac{1}{q})}\|f\|_{p}\quad
\mbox{if}\quad{\rm supp}\hat f\subset \lambda{\cal B},
\label{2.2}\\
&\|f\|_{p}\le
C\sup_{|\beta|=|\gamma|}\lambda^{-|\gamma|}\|\partial^\beta
f\|_p\le C\|f\|_p\quad\mbox{if}\quad {\rm supp}\hat f\subset
\lambda{\cal C}.\label{2.3}
\end{align}
\end{Lemma}
{\it Proof: }\,The proof can be found in \cite{Ch1}.

\section{A new Bernstein's inequality}
\quad\, Firstly, we will give certain kind of Bernstein's inequality which can be found in
[\cite{LR}, Chapter 29].
\begin{Proposition}\label{Prop3.1}
Let $2<p<\infty$. Then there exist two positive constants $c_p$ and $C_p$ such that
for every $f\in\cal{S}'$ and every $j\in\Z$, we have
\begin{align}\label{3.1}
c_p2^{\frac{2j}{p}}\|\Delta_jf\|_p\le\|\nabla(|\Delta_jf|^{\frac{p}{2}})\|_2^{\frac{2}{p}}
\le C_p2^{\frac{2j}{p}}\|\Delta_jf\|_p.
\end{align}
\end{Proposition}
Naturally, we want to establish a generalization of (\ref{3.1})
for the fractional differential operator
$\Lambda^\alpha$$(0<\alpha<1)$ which is defined by $\Lambda^\alpha
f={\cal F}^{-1}(|\xi|^\alpha\widehat{f}).$ However it seems
nontrivial, since for $p>2$, the spectrum of
$|\Delta_jf|^{\frac{p}{2}}$ can't be included in a ring although
$\textrm{supp}\,\widehat{\Delta_jf}$ is localized in $|\xi|\approx
2^j$. This section is devoted to prove Theorem \ref{Thm1}. For
this purpose, we first need the following priori lemma.

\begin{Lemma}\label{Lem3.2}
Let $p\in[1,\infty), s\in[0,p)\cap[0,2)$. Suppose that $\ell, r, m$  satisfy
$$1<\ell\le r<\infty,\quad 1<m<\infty, \quad\frac{1}{\ell}=\frac{1}{r}+\frac{p-1}{m}.$$
Then for $f(u)=|u|^p$, the following estimate holds:
\begin{align}\label{3.2}
\|f(z)\|_{\dot{B}^s_{\ell,2}}\le C_p\|z\|^{p-1}_{\dot{B}^0_{m,2}}
\|z\|_{\dot{B}^s_{r,2}}.
\end{align}
\end{Lemma}
{\it Proof:}\, Let us first recall the equivalence norm of Besov spaces: for
$0\le s<2$, $1\le \ell,q\le \infty$
$$\|v\|_{\dot{B}^s_{\ell,q}}\triangleq\bigg(\int_0^\infty t^{-sq}
\sup_{|y|\le t}\|\tau_{+y}v+\tau_{-y}v-2v\|^q_{\ell}\frac{dt}{t}\bigg)^{\frac{1}{q}},$$
where $\tau_{\pm y}v(x)=v(x\pm y)$. In the special case when $0\le s<1$, we also have
$$\|v\|_{\dot{B}^s_{\ell,q}}\triangleq\bigg(\int_0^\infty t^{-sq}
\sup_{|y|\le t}\|\tau_{+y}v-v\|^q_{\ell}\frac{dt}{t}\bigg)^{\frac{1}{q}}.$$
It is not difficult to check that

\begin{equation}\label{3.3}
|f^{[s]}(z_1)-f^{[s]}(z_2)|\le C\left\{\begin{array}{l}
(|z_1|^{p-[s]-1}+|z_2|^{p-[s]-1})|z_1-z_2|, \quad  p\geq [s]+1, \\
|z_1-z_2|^{p-[s]}, \quad p<[s]+1.\end{array}\right.
\end{equation}

where $f^{[s]}(z)=D_z^{[s]}f(z)$. For simplicity we
set $u_{\pm}\triangleq\tau_{\pm y}u$. We divide the proof of Lemma \ref{Lem3.2} into two cases.

{\bf Case 1}\,\, $p\ge 2$. We write
\begin{align}\label{3.4}
&\tau_yf(u)+\tau_{-y}f(u)-2f(u)=f(u_+)+f(u_-)-2f(u)\nonumber\\
&=f'(u)(u_++u_--2u)+\sum_{\pm}(u_{\pm}-u)
\int_0^1[f'(\lambda u_{\pm}+(1-\lambda)u)-f'(u)]d\lambda,
\end{align}
which together with (\ref{3.3}) gives that
\begin{align}
&|f(u_+)+f(u_-)-2f(u)|\nonumber\\
&\le f'(u)|u_++u_--2u|+C\sum_{\pm}|u_{\pm}-u|^2
\{\max(|u_{\pm}|, |u|)\}^{p-2}\nonumber.\end{align}
Using the H\"{o}lder inequality, we have
\begin{align}
&\|f(u_+)+f(u_-)-2f(u)\|_{\ell}\nonumber\\
&\le \|u\|_m^{p-1}\|u_++u_--2u\|_r+C\sum_{\pm}\|u_{\pm}-u\|^2_{2\theta}
\|u\|_m^{p-2},\nonumber
\end{align}
where $\theta=\frac{mr}{m+r}$.
Then by the previous equivalence norm of Besov spaces, we have
$$
\|f(u)\|_{\dot{B}^s_{\ell,2}}\le C\|u\|_{\dot{B}^s_{r,2}}\|u\|_m^{p-1}
+\|u\|_m^{p-2}\|u\|^2_{\dot{B}^{\frac{s}{2}}_{2\theta,4}}.
$$
Thanks to the interpolation inequality
$$\|u\|^2_{\dot{B}^{\frac{s}{2}}_{2\theta,4}}\le
\|u\|_{\dot{B}^s_{r,2}}\|u\|_{\dot{B}^0_{m,\infty}},$$
and the inclusion map $L^m\hookrightarrow \dot{B}^0_{m,\infty},$
we obtain
\beq\label{3.5}
\|f(u)\|_{\dot{B}^s_{\ell,2}}\le C\|u\|_{\dot{B}^s_{r,2}}\|u\|_m^{p-1}.
\eeq

{\bf Case 2}\,\, $p\le2.$  (\ref{3.3}) and (\ref{3.4}) imply that
\beq |f(u_+)+f(u_-)-2f(u)| \le
f'(u)|u_++u_--2u|+C\sum_{\pm}|u_{\pm}-u|^p.\nonumber \eeq In the
same way as leading to  (\ref{3.5}), we can deduce that
\begin{align}\label{3.6}
&\|f(u)\|_{\dot{B}^s_{\ell,2}}\le C(\|u\|^{p-1}_{m}\|u\|_{\dot{B}^s_{r,2}}
+\|u\|^p_{\dot{B}^{\frac{s}{p}}_{\ell p,2p}})\nonumber\\
&\le C(\|u\|^{p-1}_{m}\|u\|_{\dot{B}^s_{r,2}}
+\|u\|_{\dot{B}^s_{r,2}}\|u\|^{p-1}_{\dot{B}^0_{m,\infty}})
\le C\|u\|^{p-1}_{m}\|u\|_{\dot{B}^s_{r,2}}.
\end{align}
Collecting (\ref{3.5}) and (\ref{3.6}), the lemma is proved.\endproof

\begin{Remark}\label{Rem3.1} In fact,
the inequality holds for all $p\in[1,\infty), s\in[0,p)$.
But in order to make the presentation lighter, we only give the proof of the case
$s\in[0,p)\cap[0,2)$, and the  other cases can be treated in the same way.
\end{Remark}

Now let's come back to the proof of Theorem \ref{Thm1}. By homogeneity and scaling,
it is enough to prove the inequality for $j=0$.
According to the definition of Besov spaces, we have
\begin{equation}\label{3.7}\|\Lambda^\alpha(|\Delta_0f|^{\frac{p}{2}})\|_2
\cong\||\Delta_0f|^{\frac{p}{2}}\|_{\dot{B}^\alpha_{2,2}}.
\end{equation}
Applying Lemma \ref{Lem3.2} to the right hand side of (\ref{3.7})
yields that for $2\le p<\infty$, $\alpha\in[0,1]$
\begin{equation}\label{3.8}
\||\Delta_0f|^{\frac{p}{2}}\|_{\dot{B}^\alpha_{2,2}}\le
C_p\|\Delta_0f\|_{\dot{B}^0_{p,2}}^{\frac{p}{2}-1}\|\Delta_0f\|_{\dot{B}^\alpha_{p,2}}.
\end{equation}
Since ${\rm supp}\,\widehat{\Delta_0f}$ is localized in $\cal C$, by Lemma \ref{Lem2.2},
we infer that
\begin{equation}\label{3.9}
\|\Delta_0f\|_{\dot{B}^0_{p,2}},\quad  \|\Delta_0f\|_{\dot{B}^\alpha_{p,2}}
\le C\|\Delta_0f\|_{p}.\quad
\end{equation}
Collecting (\ref{3.7})-(\ref{3.9}) implies that
\begin{equation}\label{3.10}
\|\Lambda^\alpha(|\Delta_0f|^{\frac{p}{2}})\|_2^{\frac{2}{p}}
\le C_p\|\Delta_0f\|_p.
\end{equation}
In order to prove the inverse inequality, we first use Proposition \ref{Prop3.1} to get
\beq\label{3.11}
c_p\|f_0\|_p^\frac p 2\le\|\Lambda(|f_0|^{\frac{p}{2}})\|_2,
\eeq
where $f_{0}\triangleq \Delta_0f.$ To estimate $\|\Lambda(|f_0|^\frac{p}{2})\|_2$,
we decompose $\Lambda(|f_0|^\frac{p}{2})$ into
$$\Lambda(|f_{0}|^{\frac{p}{2}})=\sum_{k\ge M}\Lambda\Delta_k(|f_0|^\frac{p}{2})
+\Lambda\sum_{k<M}\Delta_k(|f_0|^\frac{p}{2})\triangleq \Lambda P_{\ge M}(|f_0|^\frac{p}{2})
+\Lambda P_{< M}(|f_0|^\frac{p}{2}),$$
for a sufficiently large $M$ which will be determined later.
 On the one hand,
we write
$$
\|\Lambda P_{\ge M}(|f_0|^\frac{p}{2})\|_2=\|\Lambda^{-\varepsilon}\Lambda^{1+\varepsilon}
(P_{\ge M}|f_0|^\frac{p}{2})\|_2,
$$
for a small enough $\varepsilon>0$ such that $1+\varepsilon<\frac p 2$.
Thanks to Lemma \ref{Lem2.2}, we get
$$
\|\Lambda^{-\varepsilon}\Lambda^{1+\varepsilon}
(P_{\ge M}|f_0|^\frac{p}{2})\|_2\le C_p2^{-M\varepsilon}\|\Lambda^{1+\varepsilon}
(|f_0|^\frac{p}{2})\|_2\approx C_p2^{-M\varepsilon}\||f_0|^\frac{p}{2}\|_{\dot B^{1+\varepsilon}_{2,2}},
$$
which together with Lemma \ref{Lem3.2} implies that
\begin{align}\label{3.12}
\|\Lambda P_{\ge M}(|f_0|^\frac{p}{2})\|_2\le C_p2^{-M\varepsilon}\|f_0\|^{\frac{p}{2}}_p.
\end{align}
On the other hand, using Lemma \ref{Lem2.2} again, we obtain
\begin{align}\label{3.13}
&\|\Lambda P_{<M}(|f_0|^\frac{p}{2})\|_2=\|\Lambda^{1-\alpha}\Lambda^{\alpha}
(P_{<M}|f_0|^\frac{p}{2})\|_2\nonumber\\ &\le C_p2^{M(1-\alpha)}
\|\Lambda^\alpha(|f_0|^{\frac{p}{2}})\|_2,
\end{align}
Combining (\ref{3.11})-(\ref{3.13}) yields that
\begin{align}
&c_p\|f_0\|^{\frac{p}{2}}_p\le\|\Lambda(|f_0|^{\frac{p}{2}})\|_2\le
\|\Lambda P_{\ge M}(|f_0|^\frac{p}{2})\|_2+\|\Lambda P_{<M}(|f_0|^\frac{p}{2})\|_2\nonumber\\
&\le
C_p\bigg(2^{-M\varepsilon}\|f_0\|^{\frac{p}{2}}_p+2^{M(1-\alpha)}
\|\Lambda^\alpha(|f_0|^{\frac{p}{2}})\|_2\bigg)\nonumber.
\end{align}
If we choose $M$ such that $C_p2^{-M\varepsilon}\le \frac{1}{2}c_p$, we conclude that
\begin{equation}\label{3.14}
c_p\|f_0\|^{\frac{p}{2}}_p\le\|\Lambda^\alpha(|f_0|^{\frac{p}{2}})\|_2.
\end{equation}
This completes the proof of Theorem \ref{Thm1}. \endproof

Finally let us recall the following improved positivity Lemma.
\begin{Lemma}\label{Lem3.3}
Suppose that $s\in[0,2]$, and $f, \Lambda^sf\in L^p(\R^2)$,
$p\ge2$. Then \beq\label{3.15} \int_{\R^2}|f|^{p-2}f\Lambda^s
fdx\ge
\frac{2}{p}\int_{\R^2}(\Lambda^{\frac{s}{2}}|f|^{\frac{p}{2}})^2dx.
\eeq
\end{Lemma}
{\it Proof:}\, The proof can be found in \cite{Ju2}.

\section{The proof of Theorem \ref{Thm2}}

In this section, we will prove Theorem \ref{Thm2}. We divided it into two parts.

\subsection{Global well-posedness for small initial data}

\quad\,\;{\bf Step 1. A priori estimates}\\

Taking the operator $\Delta_j$ on both sides of (\ref{1.1}), we have
\beq
\pa_t\Delta_j\theta+\kappa\Lambda^{2\alpha}\Delta_j\theta+u\cdot\nabla\Delta_j\theta=
[u, \Delta_j]\cdot\nabla\theta\nonumber.
\eeq
Multiplying by $p|\Delta_j\theta|^{p-2}\Delta_j\theta$ and integrating with respect to $x
$ yield that
\begin{align}\label{4.1}
&\frac{d}{dt}\|\Delta_j\theta\|^p_p+\kappa p\int_{\R^2}\Lambda^{2\alpha}\Delta_j\theta
|\Delta_j\theta|^{p-2}\Delta_j\theta dx+p\int_{\R^2} u\cdot\nabla\Delta_j\theta
|\Delta_j\theta|^{p-2}\Delta_j\theta dx\nonumber\\&=p\int_{\R^2}[u, \Delta_j]\cdot\nabla\theta
|\Delta_j\theta|^{p-2}\Delta_j\theta dx.
\end{align}
Since $\textrm{div}u=0$, by integration by parts we infer that
\beq\label{4.2}\int_{\R^2} u\cdot\nabla\Delta_j\theta
|\Delta_j\theta|^{p-2}\Delta_j\theta dx=0.\eeq
Thanks to Lemma \ref{Lem3.3} and Theorem \ref{Thm1},
we deduce that
\begin{align}\label{4.3}
p\int_{\R^2}\Lambda^{2\alpha}\Delta_j\theta
|\Delta_j\theta|^{p-2}\Delta_j\theta dx\ge2\int_{\R^2}\big(\Lambda^{\al}
|\Delta_j\theta |^{\frac{p}{2}}\big)^2dx\ge c_p2^{2\al j}\|\Delta_j\theta \|_p^p.
\end{align}
Summing up (\ref{4.1})--(\ref{4.3}) and H\"{o}lder inequality yield that
\begin{align}
\frac{d}{dt}\|\Delta_j\theta\|_p+2\kappa c_p2^{2\al j}\|\Delta_j\theta \|_p\le
C\|[u, \Delta_j]\cdot\nabla\theta\|_p.\nonumber
\end{align}
which together with Gronwall's inequality implies that
\beq\label{4.4}
\|\Delta_j\theta\|_p\le e^{-\kappa c_p t2^{2\al j}}\|\Delta_j\theta_0\|_p+
Ce^{-\kappa c_p t2^{2\al j}}\ast\|[u, \Delta_j]\cdot\nabla\theta\|_p.
\eeq
where the sign $\ast$ denotes the convolution of functions defined in $\R^+$, in details
$$e^{-\kappa c_p t2^{2\al j}}\ast f\triangleq \int_0^t
e^{-\kappa c_p (t-\tau)2^{2\al j}}f(\tau)d\tau.$$
Taking the $L^r(0,T)$  norm, $1\le r\le\infty$, $T\in(0,\infty]$, and using Young's inequality
to obtain\begin{align}\label{4.5}
\|\Delta_j\theta\|_{L^r_T(L^p)}\le\|e^{-\kappa c_p t2^{2\al j}}\|_{L^r_T}
\big(\|\Delta_j\theta_0\|_p+
C\big\|[u, \Delta_j]\cdot\nabla\theta\big\|_{L^1_T(L^p)}\big).
\end{align}
Multiplying $2^{j\sigma}$ on both sides of (\ref{4.5}), then taking $\ell^q(\Z)$ norm, we obtain
\begin{align}\label{4.6}
\|\theta\|_{\widetilde{L}^r(\dot{B}^{\sigma+\frac{2\al}{r}}_{p,q})}\lesssim
\kappa^{-1/r}\big(\|\theta_0\|_{\dot{B}^\sigma_{p,q}}+\big\|2^{j\sigma}
\|[u, \Delta_j]\cdot\nabla\theta\|_{L^1(\R^+,
L^p)}\big\|_{\ell^q(\Z)}\big),
\end{align}
where we used the fact that \beq\label{4.7} \big\|e^{-\kappa c_p
t2^{2\al j}}\big\|_{L^{r}_T}\le\bigg( \frac{1-e^{-r\kappa c_p
2^{2\al j}T}}{r\kappa c_p2^{2\al j}}\bigg)^\frac{1}{r},
\qquad{\hbox{for}}\quad 1\le r\le\infty, \eeq and
$\sigma=\frac{2}{p}+1-2\al.$ On the other hand, it follows from
Proposition \ref{PropA.3} that
\begin{align}\label{4.8}
\big\|2^{j\sigma} \|[u, \Delta_j]\cdot\nabla\theta\|_{L^1(\R^+,
L^p)}\big\|_{\ell^q(\Z)} &\le C
\|u\|_{\widetilde{L}^2(\dot{B}^{\frac{2}{p}+1-\al}_{p,q})}
\|\theta\|_{\widetilde{L}^2(\dot{B}^{\frac{2}{p}+1-\al}_{p,q})}\nonumber\\
&\le C \|\theta\|_{\widetilde{L}^\infty(\dot B^\sigma_{p,q})}
\|\theta\|_{\widetilde{L}^1(\dot{B}^{\frac{2}{p}+1}_{p,q})},
\end{align}
where in the last inequality we have used the interpolation and
the fact that \beq\label{4.9} \|u\|_{L^r(\dot{B}^s_{p,q})}=\|{\cal
R}_k\theta\|_{L^r(\dot{B}^s_{p,q})} \le
C\|\theta\|_{L^r(\dot{B}^s_{p,q})},\quad \textrm{for}\,\, s\in
\R,\,\, (r,p,q)\in [1,\infty]^3, \eeq since $\|\Delta_j{\cal
R}_k\theta\|_{p}\thickapprox\|\widetilde{\Delta}_{j}{\cal R}_k\Delta_j\theta\|_{p}\le
C\|\Delta_j\theta\|_{p}$  for all $1\le p\le \infty$, here
$\widetilde{\Delta}_{j}=(\Delta_{j-1}+\Delta_{j}+\Delta_{j+1})$.
Combining (\ref{4.6}) and (\ref{4.8}), we get
\begin{align}\label{4.10}
\|\theta\|_{\widetilde{L}^r(\dot{B}^{\sigma+\frac{2\al}{r}}_{p,q})}\lesssim
\kappa^{-1/r}\big(\|\theta_0\|_{\dot{B}^\sigma_{p,q}}
+C \|\theta\|_{\widetilde{L}^\infty(\dot B^\sigma_{p,q})}
\|\theta\|_{\widetilde{L}^1(\dot{B}^{\frac{2}{p}+1}_{p,q})}\big).
\end{align}
On the other hand, it follows from  (\cite{Cor3}, Corollary 2.6) that
\beq\label{4.11}
\|\theta(t,x)\|_p\le \|\theta_0(x)\|_p,\quad t\ge 0, \eeq
which together with (\ref{4.10}) implies that
\begin{align}\label{4.12}
&\|\theta(t)\|_{\widetilde{L}^\infty(
B^\sigma_{p,q})}+c_1\kappa
\|\theta(t)\|_{\widetilde{L}^1(\dot B^{\frac{2}{p}+1}_{p,q})}
\le2\|\theta_0\|_{{B}^\sigma_{p,q}}+C\|\theta\|_{\widetilde{L}^\infty({B}^{\sigma}_{p,q})}
\|\theta\|_{\widetilde{L}^1(\dot{B}^{\frac{2}{p}+1}_{p,q})}.
\end{align}

{\bf Step 2. Approximation solutions and uniform estimates}\\

Let us define the sequence $\{\th^{(n)}, u^{(n)}\}_{n\in\N_0}$ by
the following systems:
\begin{equation}\label{4.13}
\left\{
\begin{aligned}
&\partial_t\theta^{(n+1)}+u^{(n)}\cdot \nabla
\theta^{(n+1)}+\kappa(-\Delta)^{\alpha} \theta^{(n+1)}=0,\quad
x\in \R^2,\, t>0,\\& u^{(n)}={\cal R}^\bot\theta^{(n)},
\\
&\theta^{(n+1)}(0,x)=\theta_0^{(n+1)}(x)=\sum_{j\le
n+1}\Delta_j\theta_0(x).
\end{aligned}
\right.
\end{equation}
Setting $(\th^{(0)}, u^{(0)})$=$(0,0)$ and solving the linear
system, we can find $\{\th^{(n)}, u^{(n)}\}_{n\in\N_0}$ for all
$n\in \N_0$. As in Step 1, we can deduce that
\begin{align}\label{4.14}
&\|\theta^{(n+1)}(t)\|_{\widetilde{L}^\infty(\dot
B^\sigma_{p,q})}+c_1\kappa
\|\theta^{(n+1)}(t)\|_{\widetilde{L}^1(\dot B^{\frac{2}{p}+1}_{p,q})}\nonumber\\
&\le2\|\theta_0^{(n+1)}\|_{\dot
B^\sigma_{p,q}}+C_2{(c_1\kappa)}^{-1}
\big(\|\theta^{(n)}\|_{\widetilde{L}^\infty(\dot
B^{\sigma}_{p,q})}+c_1\kappa
\|\theta^{(n)}\|_{\widetilde{L}^1\big(\dot
B^{\frac{2}{p}+1}_{p,q})}\big)\nonumber\\& \quad\,\times\big(
\|\theta^{(n+1)}\|_{\widetilde{L}^\infty(\dot
B^{\sigma}_{p,q})}+c_1\kappa
\|\theta^{(n+1)}\|_{\widetilde{L}^1(\dot
B^{\frac{2}{p}+1}_{p,q})}\big).
\end{align}
If we take $\epsilon>0$ such that $\|\theta_0\|_{\dot
B^\sigma_{p,q}}\le\epsilon\kappa, \epsilon\le \frac{c_1} {8C_2}$,
then for all $n$, we will show \beq\label{4.15}
\|\theta^{(n)}(t)\|_{\widetilde{L}^\infty(\dot
B^\sigma_{p,q})}+c_1\kappa
\|\theta^{(n)}(t)\|_{\widetilde{L}^1(\dot
B^{\frac{2}{p}+1}_{p,q})} \le 4\|\theta_0\|_{\dot B^\sigma_{p,q}}.
\eeq In fact, assume that
$\|\theta^{(k)}\|_{\widetilde{L}^\infty(\dot
B^\sigma_{p,q})}+c_1\kappa \|\theta^{(k)}\|_{\widetilde{L}^1(\dot
B^{\frac{2}{p}+1}_{p,q})}\le 4\|\theta_0\|_{\dot B^\sigma_{p,q}}$
for $k=0,\cdots,n$. It follows from (\ref{4.14}) that
\begin{align}\label{4.16}
&\|\theta^{(n+1)}\|_{\widetilde{L}^\infty(\dot
B^\sigma_{p,q})}+c_1\kappa
\|\theta^{(n+1)}\|_{\widetilde{L}^1(\dot B^{\frac{2}{p}+1}_{p,q})}
\nonumber\\&\le2\|\theta_0\|_{\dot
B^\sigma_{p,q}}+C_2(c_1\kappa)^{-1}4 \|\theta_0\|_{\dot
B^\sigma_{p,q}} \big(\|\theta^{(n+1)}\|_{\widetilde{L}^\infty(\dot
B^{\sigma}_{p,q})}+c_1\kappa
\|\theta^{(n+1)}\|_{\widetilde{L}^1(\dot
B^{\frac{2}{p}+1}_{p,q})}\big) \nonumber\\&\le2\|\theta_0\|_{\dot
B^\sigma_{p,q}}+\frac{1}{2}
\big(\|\theta^{(n+1)}\|_{\widetilde{L}^\infty(\dot
B^{\sigma}_{p,q})}+c_1\kappa
\|\theta^{(n+1)}\|_{\widetilde{L}^1(\dot
B^{\frac{2}{p}+1}_{p,q})}\big),
\end{align}
which implies (\ref{4.15}). Summing up (\ref{4.11}) and
(\ref{4.15}), we finally get for all $n$, \beq\label{4.17}
\|\theta^{(n)}(t)\|_{\widetilde{L}^\infty(
B^\sigma_{p,q})}+c_1\kappa
\|\theta^{(n)}(t)\|_{\widetilde{L}^1(\dot
B^{\frac{2}{p}+1}_{p,q})} \le4\|\theta_0\|_{ B^\sigma_{p,q}}. \eeq

{\bf Step 3. Compactness arguments and Existence}\\

We will show that, up to a subsequence, the sequence
$\{\theta^{(n)}\}$ converges in ${\cal D}'(\R^+\times \R^2)$ to a
solution $\theta$ of (\ref{1.1}). The proof is based on
compactness arguments. First we show that $\partial_t\theta^{(n)}$
is uniformly bounded in the  space ${L}^\infty(B^{-2\al}_{p,q})$.
By (\ref{4.13}), $\partial_t\theta^{(n+1)}$ satisfies the equation
$$\partial_t\theta^{(n+1)}=-\nabla \cdot(u^{(n)}
\theta^{(n+1)})-\kappa(-\Delta)^\al \theta^{(n+1)}.$$ Then thanks
to Proposition \ref{PropA.1} with $p\not=\infty$, we get
\begin{align}
\|\partial_t\theta^{(n+1)}\|_{{L}^\infty(B^{-2\al}_{p,q})}&\lesssim
\|\theta^{(n+1)}\|_{{L}^\infty(B^{0}_{p,q})}+
\|u^{(n)}\|_{{L}^\infty(L^p)}
\|\theta^{(n+1)}\|_{{L}^\infty({B}^{\sigma}_{p,q})}\nonumber\\&\quad+
\|\theta^{(n+1)}\|_{{L}^\infty(L^p)}\|u^{(n)}\|_{{L}^\infty({B}^{\sigma}_{p,q})}\nonumber\\
&\lesssim
\|\theta^{(n+1)}\|_{{L}^\infty(B^{\sigma}_{p,q})}+\|\theta^{(n)}\|_{{L}^\infty({B}^{\sigma}_{p,q})}
\|\theta^{(n+1)}\|_{{L}^\infty({B}^{\sigma}_{p,q})}
<\infty,\nonumber
\end{align}
where we have used the fact: for $s>0$,
${B}^s_{p,q}=L^p\cap\dot{B}^s_{p,q}$, and the inclusion map
$B^\sigma_{p,q}\subset B^0_{p,q}$. We remark that the above
inequality can be obtained also by Proposition \ref{PropA.2} with
$s=-2\sigma$, $s_1$ be an  any number such that $0<s_1<\frac 2p$.
Now let us turn to the proof of the existence. Observe that for
any $\chi\in C_c^\infty(\R^2)$, the map: $u\mapsto \chi u$ is
compact from $B^\sigma_{p,q}(\R^2)$ into $L^p(\R^2)$. This can be
proved by noting that the map: $u\mapsto \chi u$ is compact from
$H^{s'}_p$ into $H^{s}_p$ for $s'>s$, $p<\infty$, and the
embedding relation $B^\sigma_{p,q}\hookrightarrow
B^{\sigma-\epsilon}_{p,2}\hookrightarrow H^{\sigma-\ep}_p$(by
Proposition \ref{Prop2.1}(iii)). Thus by the Lions-Aubin
compactness theorem(see \cite{Tem}), we can conclude that there
exists a subsequence $\{\theta^{(n_k)}\}$ and a function $\theta $
so that
$$
\lim_{n_k\rightarrow+\infty}\theta^{(n_k)}= \theta \quad
\textrm{in} \quad L^p_{loc}(\R^+\times \R^2).
$$
Moreover,  the
uniform estimate (\ref{4.17}) allows us to
conclude that
$$\theta(t,x)\in \widetilde{L}^\infty(0,\infty; B^\sigma_{p,q})\cap
 \widetilde{L}^1(0,\infty;\dot B^{\frac{2}{p}+1}_{p,q}),$$
and
$$
\|\theta(t)\|_{\widetilde{L}^\infty
(B^\sigma_{p,q})}+\|\theta(t)\|_{\widetilde{L}^1(\dot B^{\frac{2}{p}+1}_{p,q})}\le
4\|\theta_0\|_{B^\sigma_{p,q}}.
$$
Then by a standard limit argument, we can prove that the limit
function $\theta(t,x)$ satisfies the equation (\ref{1.1}) in the
sense of distribution.

We still have to prove $\theta(t,x)$
belongs to  $C(\R^+; B^\sigma_{p,q})$. Our idea comes from \cite{Dan2}. We observe that
\beq\label{4.18}
\pa_t\Delta_j\theta=-\kappa\Lambda^{2\al}\Delta_j\theta-\Delta_j\na\cdot(u\theta)
\eeq
For fixed $j$, the right hand side of (\ref{4.18}) belongs to
$L^\infty(0,\infty; B^\sigma_{p,q})$, which can be easily proved by using Lemma 2.2.
Therefore, we infer that $\pa_t\Delta_j\theta\in L^\infty(0,\infty; B^\sigma_{p,q})$
for fixed $j$,
which implies that each $\Delta_j\theta$ is continuous in time in  $B^\sigma_{p,q}$.
On the other hand, note that $$\|\theta\|_{\widetilde{L}^\infty(B^\sigma_{p,q})}=\bigg(
\sum_{j\in\Z}\sup_{t\ge0}\big(2^{j\sigma}\|\Delta_j\theta\|_{L^p}\big)^q\bigg)^{\frac{1}{q}}
<+\infty,$$
which implies that $\sum_{|j|\le n}\Delta_j\theta $ converges uniformly in
$L^\infty(\R^+; B^\sigma_{p,q})$ to $\theta(t,x)$. Hence,
$\theta(t,x)\in C(\R^+; B^\sigma_{p,q}).$\\

{\bf Step 4. Uniqueness}\\

Assume that $\theta'\in \widetilde L^\infty( B^\sigma_{p,q})\cap
\widetilde L^1(\dot B^{\frac{2}{p}+1}_{p,q})$ is another solution of (\ref{1.1}) with the same
initial data $\theta_0(x)$. Let $\delta\theta=\theta-\theta'$ and
$\delta u=u-u'$. Then $(\delta \theta,\delta u)$ satisfy the
following equations
\begin{equation}\label{4.20} \left\{
\begin{aligned}
&\partial_t\delta\theta+u\cdot \nabla \delta\theta+\delta u\cdot
\nabla\theta'+\kappa(-\Delta)^\alpha
\delta\theta=0,\quad x\in \R^2,\, t>0,\\
&\delta u={\cal R}^\bot\delta\theta,
\\ &\delta\theta(0,x)=0.
\end{aligned}
\right.
\end{equation}
Following the same way as a priori estimates, we can deduce that
$$
\frac{d}{dt}\|\Delta_j\delta\theta\|_p+2\kappa c_p2^{2\al j}\|\Delta_j\delta\theta \|_p\le
C\big(\|[u, \Delta_j]\cdot\nabla\delta\theta\|_p+\|\Delta_j(\delta u\cdot \nabla\theta')\|_p\big),
$$
which together with Gronwall's inequality leads to
\begin{align}\label{4.19}
\|\Delta_j\delta\theta\|_p\le
Ce^{-\kappa c_p t2^{2\al j}}\ast\big(\|[u, \Delta_j]\cdot\nabla\delta\theta\|_p
+\|\Delta_j(\delta u\cdot \nabla\theta')\|_p\big).
\end{align}
Choose a positive number $\eta$ such that $\frac{2\alpha}p<\eta<\frac{2}{p}$.
Thanks to Proposition \ref{PropA.2}, (\ref{4.9}), and  interpolation, we get
\begin{align}\label{4.20}
&\|\delta u\cdot\nabla \theta'\|_{\widetilde{L}^1_T(\dot{B}^{\frac{2}{p}-\eta}_{p,q})}
\lesssim
\|\delta u\|_{\widetilde{L}^{p}_T(\dot{B}^{\frac{2}{p}-\eta+\frac{2\alpha}p}_{p,q})}
\|\theta'\|_{\widetilde{L}^{p'}_T(\dot{B}^{\frac{2}{p}+1-\frac{2\alpha}p}_{p,q})}\nonumber\\&\lesssim
\|\delta \theta\|_{\widetilde{L}^{p}_T(\dot{B}^{\frac{2}{p}-\eta+\frac{2\alpha}p}_{p,q})}
\|\theta'\|_{\widetilde{L}^{p'}_T(\dot{B}^{\frac{2}{p}+1-\frac{2\alpha}p}_{p,q})}\nonumber\\&\lesssim
\|\delta\theta\|_{\widetilde{L}^\infty_T(\dot{B}^{\frac{2}{p}-\eta}_{p,q})}^{\frac{1}{p'}}\|\delta\theta\|
_{\widetilde{L}^1_T(\dot{B}^{\frac{2}{p}-\eta+2\alpha}_{p,q})}^{\frac{1}{p}}
\|\theta'\|_{\widetilde{L}^{p'}_T(\dot{B}^{\frac{2}{p}+1-\frac{2\alpha}p}_{p,q})}.
\end{align}
On the other hand, thanks to Proposition \ref{PropA.3}, (\ref{4.9}) we have
\begin{align}\label{4.21}
&\|[u, \Delta_j]\cdot\nabla\delta\theta\|_{L^1_T(L^p)}\lesssim c_j2^{-j(\frac{2}{p}-\eta)}
\|\theta\|_{\widetilde{L}^{2}_T(\dot{B}^{\frac{2}{p}+1-\al}_{p,q})}
\|\delta \theta\|_{\widetilde{L}^{2}_T(\dot{B}^{\frac{2}{p}-\eta+\al}_{p,q})}\nonumber\\
&\lesssim c_j2^{-j(\frac{2}{p}-\eta)}
\|u\|_{\widetilde{L}^{2}_T(\dot{B}^{\frac{2}{p}+1-\al}_{p,q})}
\|\delta \theta\|_{\widetilde{L}^{\infty}_T(\dot{B}^{\frac{2}{p}-\eta}_{p,q})}^{\frac{1}{2}}
\|\delta \theta\|_{\widetilde{L}^{1}_T(\dot{B}^{\frac{2}{p}-\eta+2\al}_{p,q})}^{\frac{1}{2}}.
\end{align}
where $\|c_j\|_{\ell^q(\Z)}\le1$.  Taking $L^\infty$($L^1$, respectively) norm on time,
and using Young's inequality, then
multiplying $2^{j(\frac{2}{p}-\eta)}$($2^{j(\frac{2}{p}-\eta+2\al)}$, respectively)
 on both sides of (\ref{4.19}),
then taking $\ell^q(\Z)$ norm, we have
\begin{align}\label{4.22}
&Z(T)\triangleq\|\delta\th\|_{\widetilde{L}^\infty_T(\dot{B}^{\frac{2}{p}-\eta}_{p,q})}
+\|\delta\theta\|
_{\widetilde{L}^1_T(\dot{B}^{\frac{2}{p}-\eta+2\al}_{p,q})}\nonumber\\&
\lesssim\big\|2^{j(\frac{2}{p}-\eta)}
\|[u, \Delta_j]\nabla\delta\theta\|_{L^1_T(L^p)}\big\|_{\ell^q(\Z)}+
\|\delta u\cdot \nabla\theta'\|_{\widetilde{L}^{1}_T(\dot{B}^{\frac{2}{p}-\eta}_{p,q})}
\nonumber\\
&\lesssim\bigg(\|\theta'\|_{\widetilde{L}^{p'}_T(\dot{B}^{\frac{2}{p}+1-\frac{2\alpha}p}_{p,q})}
+\|\theta\|_{\widetilde{L}^{2}_T(\dot{B}^{\frac{2}{p}+1-\al}_{p,q})}\bigg)
Z(T),
\end{align}
where we have used (\ref{4.20}) and (\ref{4.21}) in the last inequality.
Now it is clear that two terms in the bracket of the right hand side of (\ref{4.22}) tend to
0 as $T$ goes to 0. Therefore, if $T$ has been chosen small enough, then it follows from (\ref{4.22}) that
$Z\equiv0$ on $[0,T]$ which implies that $\delta\theta\equiv0$. Then by a standard
continuous argument, we can show that $\delta\theta(t,x)=0$ in
$[0,+\infty)\times \R^2$, i.e. $\theta(t,x)=\theta'(t,x)$.

This completes the proof of global well-posedness.\endproof

\subsection{Local well-posedness for large initial data}
\quad\, Now we prove the local well-posedness for the large initial data. As the existence result will be essentially followed
from the a priori estimate. For simplicity, we only
present a priori estimate of the solution $\theta(t,x)$.

Let us return to (\ref{4.5}). Taking $r=2$ in (\ref{4.5}),
multiplying $2^{j(\frac{2}{p}+1-\al)}$ on both sides of
(\ref{4.5}), then taking $\ell^q(\Z)$ norm and  applying
Proposition \ref{PropA.3} and (\ref{4.9}), we get
\begin{align}\label{4.23}
\|\theta\|_{\widetilde{L}^2_T(\dot{B}^{\frac{2}{p}+1-\al}_{p,q})}&\le
C_3\kappa^{-\frac 1 2}\bigg(\big\|E_{j}(T)^{\frac12}2^{j\sigma}\|\Delta_j\theta_0\|_{p}\big\|_{\ell^q(\Z)}+
\|u\|_{\widetilde{L}^2_T(\dot{B}^{\frac{2}{p}+1-\al}_{p,q})}
\|\theta\|_{\widetilde{L}^2_T(\dot{B}^{\frac{2}{p}+1-\al}_{p,q})}\bigg)\nonumber\\ &\le C_3
\kappa^{-\frac 1 2}\bigg(\big\|E_{j}(T)^{\frac12}2^{j\sigma}\|\Delta_j\theta_0\|_{p}\big\|_{\ell^q(\Z)}+
\|\theta\|_{\widetilde{L}^2_T(\dot{B}^{\frac{2}{p}+1-\al}_{p,q})}
^2\bigg),
\end{align}
where $$E_{j}(T)\triangleq 1-e^{-\kappa c_p 2^{2\al j}T}.$$
Set
$$
T_0\triangleq\sup\bigg\{T'>0;\,\bigg(\sum_{j\in\Z}E_{j}(T')^{\frac q 2}2^{j\sigma q}
\|\Delta_j\theta_0\|^q_{p}\bigg)^\frac{1}{q}\le\frac{\kappa}{2C_3^2}\bigg\}.
$$
Then the inequality (\ref{4.23}) implies that there holds for $T\in [0,T_0]$
\beno
\|\theta(t)\|_{\widetilde{L}^2_{T}(\dot{B}^{\frac{2}{p}+1-\al}_{p,q})}
\le2\big\|E_{j}(T)^{\frac12}2^{j\sigma}\|\Delta_j\theta_0\|_{p}\big\|_{\ell^q(\Z)},
\eeno
which together with (\ref{4.6}) and (\ref{4.8}) leads to
\beno
\|\theta(t)\|_{\widetilde{L}^\infty_T(\dot
B^\sigma_{p,q})}+c_1\kappa
\|\theta(t)\|_{\widetilde{L}^1_T(\dot B^{\frac{2}{p}+1}_{p,q})}
\le 2\|\theta_0\|_{\dot B^\sigma_{p,q}}+
C \|\theta\|_{\widetilde{L}^2_T(\dot B^{\frac{2}{p}+1-\al}_{p,q})}^2\le
C\|\theta_0\|_{\dot B^\sigma_{p,q}}.
\eeno
Combining with (\ref{4.11}), we obtain for $T\in [0,T_0]$
\beno
\|\theta(t)\|_{\widetilde{L}^\infty_T(
B^\sigma_{p,q})}+c_1\kappa
\|\theta(t)\|_{\widetilde{L}^1_T(\dot B^{\frac{2}{p}+1}_{p,q})}
\le C\|\theta_0\|_{B^\sigma_{p,q}}.
\eeno

This completes the proof of local well-posedness.\endproof

\section{Appendix}
Firstly, we recall the paradifferential calculus which enables us to define a generalized
product between distributions, which is continuous in many functional spaces where the usual
product does not make sense (see \cite{Bon}). The paraproduct between $u$ and $v$ is defined
by$$T_uv\triangleq\sum_{j\in\Z}S_{j-1}u\Delta_jv.$$
We then have the following formal decomposition:
\beq\label{A.1}
uv=T_uv+T_vu+R(u,v),
\eeq
with $$R(u,v)=\sum_{j\in\Z}\Delta_ju\widetilde{\Delta}_jv\quad\mbox{and}\quad
\widetilde{\Delta}_j=\Delta_{j-1}+\Delta_j+\Delta_{j+1}.$$
The decomposition (\ref{A.1}) is called the Bony's paraproduct decomposition.

Now we  state some results  about the product estimates in
Besov spaces.

\begin{Proposition}\label{PropA.1}
Let $s>-\frac{d}{p}$, $2\le p\le\infty$, $1\le q\le\infty.$ Then
\begin{align}\label{A.2}
\|uv\|_{ B^s_{p,q}}
\lesssim\|u\|_{p}\|v\|_{B^{\frac{d}{p}+s}_{p,q}}
+\|v\|_{p}\|u\|_{B^{\frac{d}{p}+s}_{p,q}}.
\end{align}
\end{Proposition}
{\it Proof: }\,  Using lemma \ref{Lem2.2}, we have
\begin{align}\label{A.3}
\|S_0(uv)\|_{p}\le \|S_0(uv)\|_{\frac{p}{2}}\le C\|u\|_{p}\|v\|_{p}\lesssim
\|u\|_{p}\|v\|_{B^{\frac{d}{p}+s}_{p,q}}.
\end{align}
Then using the  Bony's paraproduct decomposition  and the
property of quasi-orthogonality (\ref{2.1}), for fixed $j\ge0$, we have
\begin{align}\label{A.4}
\Delta_j(uv)&=\sum_{|k-j|\le4}\Delta_j(S_{k-1}u\Delta_kv)+
\sum_{|k-j|\le4}\Delta_j(S_{k-1}v\Delta_ku)
+\sum_{k\ge j-2}\Delta_j(\Delta_{k}u\widetilde{\Delta}_{k}v)\nonumber\\
&\triangleq I+II+III.\end{align}
We shall  estimate the
above three terms separately.
Using Young's inequality and  Lemma \ref{Lem2.2}, we get
\beq
\|\Delta_j(S_{k-1}u\Delta_kv)\|_{p}\lesssim
2^{\frac{d}{p}j}\|S_{k-1}u\|_{p}\|\Delta_kv\|_{p}.
\nonumber
\eeq
Thus we have
\begin{align}\label{A.5}
2^{sj}\|{I}\|_{p}\lesssim \|u\|_{p}
\sum_{|k-j|\le4}2^{(j-k)(\frac{d}{p}+s)}2^{k(\frac{d}{p}+s)}
\|\Delta_k v\|_{p} \lesssim
c_j\|u\|_{p}\|v\|_{B^{\frac{d}{p}+s}_{p,q}}.
\end{align}
where the $\ell^q(\Z)$ norm of $c_j$ is less than or equals to 1.
Similarly to $II$, we have
\ben\label{A.6}
2^{sj}\|II\|_{p}&\lesssim& c_j\|v\|_{p}\|u\|_{B^{\frac{d}{p}+s}_{p,q}}.
\een
Now we turn to estimate ${III}$. From Lemma \ref{Lem2.2},
Young's inequality, and H\"{o}lder inequality
we have \beq
\|\Delta_j(\Delta_{k}u\widetilde{\Delta}_kv)\|_{p}\lesssim
2^{j\frac{d}{p}}\|\Delta_j(\Delta_{k}u\widetilde{\Delta}_kv)\|_{\frac{p}{2}}\lesssim
2^{j\frac{d}{p}}\|\Delta_{k}u\|_{p}\|\widetilde{\Delta}_{k}v\|_{p}
.\nonumber\eeq
So, we get by $\ell^1(\Z)-\ell^q(\Z)$ convolution, \begin{align}\label{A.7}
2^{sj}\|{III}\|_{p}\lesssim \|u\|_{p}
\sum_{k\ge j-2}2^{(j-k)(\frac{d}{p}+s)}2^{k(\frac{d}{p}+s)}
\|\widetilde{\Delta}_{k}v\|_{p}
\lesssim c_j
\|u\|_{p}\|v\|_{B^{\frac{d}{p}+s}_{p,q}}.
\end{align}
where we have used the fact $s+\frac{d}{p}>0$.
Summing up (\ref{A.3}), (\ref{A.5})-(\ref{A.7}), we obtain the desired inequality
(\ref{A.2}). \endproof

\begin{Proposition}\label{PropA.2}
Let $s>-\frac{d}{p}-1$, $s<s_1<\frac{d}{p}$,
$2\le p\le\infty$, $1\le q\le\infty,$
$\frac{1}{r}=\frac{1}{r_1}+\frac{1}{r_2}=\frac{1}{r_1}+\frac{1}{r_2}\le1$,
and $u$ be a solenoidal vector field.
Then
\begin{align}\label{A.8}
\|u\cdot\na v\|_{\widetilde{L}^r_t(\dot B^s_{p,q})}
\lesssim\|u\|_{\widetilde{L}^{r_1}_t(\dot B^{s_1}_{p,q})}
\|\na v\|_{\widetilde{L}^{r_2}_t(\dot B^{s+\frac{d}{p}-s_1}_{p,q})}
.
\end{align}
If $s_1=\frac{d}{p}$ or $s_1=s$, $q$ has to be equal to 1.
\end{Proposition}
{\it Proof: }\,Throughout the proof, the summation convention over repeated indices
$i\in[1,d]$ is used. Similarly to the proof of Proposition \ref{PropA.1},
we will estimate separately each part of the Bony's paraproduct decomposition of $u^i\pa_iv$.

By Lemma \ref{Lem2.2}, we have
\begin{align}
&\|\Delta_j(S_{k-1}u^i\Delta_k\pa_i v)\|_{L^{r}_t(L^p)}\lesssim
\|S_{k-1}u\|_{L^{r_1}_t(L^\infty)}\|\Delta_k\na v\|_{L^{r_2}_t(L^p)}\nonumber\\&\lesssim
\sum_{k'\le k-2}2^{(\frac{d}{p}-s_1)k'}2^{s_1k'}\|\Delta_{k'}u\|_{L^{r_1}_t(L^p)}
\|\Delta_k\na v\|_{L^{r_2}_t(L^p)}\nonumber\\
&\lesssim2^{(\frac{d}{p}-s_1)k}\|u\|_{\widetilde{L}^{r_1}_t(\dot B^{s_1}_{p,q})}
\|\Delta_k\na v\|_{L^{r_2}_t(L^p)},
\nonumber
\end{align}
where the fact $s_1<\frac{d}{p}$ has been used in the last inequality.
Hence,  we get
\begin{align}\label{A.9}
&2^{sj}\big\|\sum_{|j-k|\le4}\Delta_j(S_{k-1}u^i\pa_i\Delta_kv)\big\|_{L^{r}_t(L^p)}
\lesssim\|u\|_{L^{r_1}_t(\dot B^{s_1}_{p,q})}\sum_{|j-k|\le4}2^{(j-k)s}2^{(s+\frac{d}{p}-s_1)k}
\|\Delta_k\na v\|_{L^{r_2}_t(L^p)}\nonumber\\&\lesssim c_j
\|u\|_{\widetilde{L}^{r_1}_t(\dot B^{s_1}_{p,q})}\|v\|_{\widetilde{L}^{r_2}_t(\dot B^{s+\frac{d}{p}+1-s_1}_{p,q})}.
\end{align}
where $\|c_j\|_{\ell^q(\Z)}\le1$.
Since ${\rm div} u=0$ and $p\ge2$, Lemma \ref{Lem2.2} applied yields that
\begin{align}
\|\Delta_j(\Delta_ku\widetilde{\Delta}_k\pa_iv)\|_{L^{r}_t(L^p)}
\lesssim2^{j(\frac{d}{p}+1)}\|\Delta_j(\Delta_ku^i\widetilde{\Delta}_kv)\|
_{L^{r}_t(L^\frac{p}{2})}.\nonumber
\end{align}
Thus by H\"{o}lder inequality and $\frac{d}{p}+1+s>0$, we have
\begin{align}\label{A.10}
&2^{sj}\big\|\sum_{k\ge j-2}
\Delta_j(\Delta_ku^i\widetilde{\Delta}_k\pa_iv)\big\|_{L^{r}_t(L^p)}
\lesssim\sum_{k\ge j-2}2^{j(\frac{d}{p}+1+s)}\|\Delta_ku\|_{L^{r_1}_t(L^p)}
\|\widetilde{\Delta}_kv\|_{L^{r_2}_t(L^p)}\nonumber\\&\lesssim
\|v\|_{\widetilde{L}^{r_2}_t(\dot B^{s+\frac{d}{p}+1-s_1}_{p,q})}
\sum_{k\ge j-2}2^{(j-k)(\frac{d}{p}+1+s)}2^{ks_1}\|\Delta_k u\|_{L^{r_1}_t(L^{p})}\nonumber\\&\lesssim c_j
\|u\|_{\widetilde{L}^{r_1}_t(\dot B^{s_1}_{p,q})}
\|v\|_{\widetilde{L}^{r_2}_t(\dot B^{s+\frac{d}{p}+1-s_1}_{p,q})}.
\end{align}
On the other hand, due to $s_1>s$, we get
\begin{align}
\|\Delta_j(\Delta_ku^iS_{k-1}\pa_i\Delta_kv)\|_{L^{r}_t(L^p)}\lesssim&
\sum_{k'\le
k-2}2^{(\frac{d}{p}+1+s-s_1)k'}2^{(s_1-s)k'}\|\Delta_{k'}v\|_{L^{r_2}_t(L^p)}
\|\Delta_ku\|_{L^{r_1}_t(L^p)}\nonumber\\&\lesssim2^{(s_1-s)k}
\|v\|_{\widetilde{L}^{r_2}_t(\dot B^{s+\frac{d}{p}+1-s_1}_{p,q})}
\|\Delta_ku\|_{L^{r_1}_t(L^p)}.\nonumber
\end{align}
Then we have
\begin{align}\label{A.11}
&2^{sj}\big\|\sum_{|j-k|\le4}\Delta_j(\Delta_ku^iS_{k-1}\pa_i\Delta_kv)\big\|_{L^{r}_t(L^p)}
\nonumber\\&\lesssim\|v\|_{L^{r_2}_t(\dot B^{s+\frac{d}{p}+1-s_1}_{p,q})}\sum_{|j-k|\le4}
2^{s(j-k)}2^{s_1k}\|\Delta_ku\|_{L^{r_1}_t(L^p)}\nonumber\\&\lesssim c_j
\|v\|_{\widetilde{L}^{r_2}_t(\dot B^{s+\frac{d}{p}+1-s_1}_{p,q})}
\|u\|_{\widetilde{L}^{r_1}_t(\dot B^{s_1}_{p,q})}.
\end{align}
Summing up (\ref{A.9})-(\ref{A.11}), the desired inequality (\ref{A.8}) is proved.\endproof

\vspace{.3cm}
Finally we give the commutator estimate.
\begin{Proposition}{\label{PropA.3}}
Let $1\le p,\, q\le\infty$, $\frac{1}{r}=\frac{1}{r_1}+\frac{1}{r_2}\le1$,
$\rho<1$, $\gamma>-1$ and $u$ be a solenoidal vector field.
Assume in addition that
$$\rho-\gamma+d\min\big(1,\frac{2}{p}\big)>0\quad\mbox{and}\quad\rho+\frac{d}{p}>0.$$
Then
the following inequality holds:
\beq\label{A.12}
\|[u,\Delta_j]\cdot\na v\|_{L^r_t(L^p)}\lesssim c_j2^{-j(\frac{d}{p}+\rho-1-\gamma)}
\|\na u\|_{\widetilde{L}^{r_1}_t(\dot B^{\frac{d}{p}+\rho-1}_{p,q})}
\|\na v\|_{\widetilde{L}^{r_2}_t(\dot B^{\frac{d}{p}-\gamma-1}_{p,q})}.
\eeq
where $c_j$ denotes a positive series with $\|c_j\|_{\ell^q(\Z)}\le1$. In the above,
we denote
$$[u,\Delta_j]\cdot\na v=\sum_{1\le i\le d}u^i\Delta_j\pa_iv-\Delta_j(u^i\pa_i v).$$
If $\rho=1$, $\|\na u\|_{\widetilde{L}^{r_1}_t(\dot B^{\frac{d}{p}+\rho-1}_{p,q})}$
has to be replaced
by $\|\na u\|_{\widetilde{L}^{r_1}_t(\dot B^{\frac{d}{p}+\rho-1}_{p,q})
\cap L^{r_1}_t(L^\infty)}.$ If
$\gamma=-1$, $\|\na v\|_{\widetilde{L}^{r_2}_t(\dot B^{\frac{d}{p}-\gamma-1}_{p,q})}$
has to be replaced
by $\|\na v\|_{\widetilde{L}^{r_2}_t(\dot B^{\frac{d}{p}-\gamma-1}_{p,q})
\cap L^{r_1}_t(L^\infty)}.$
\end{Proposition}
{\it Proof:}\, The proof is a straightforward adaptation of Lemma A.1 in \cite{Dan3} which
is  a version of the commutator estimate in  Besov space.

\textbf{Acknowledgements} We would like  to thank  Professors H.
Smith and T. Tao so much for  their  helpful discussion and
suggestions. The authors are also deeply grateful to the referee
for their valuable advices. Q. Chen and C. Miao  were partly
supported by the National Natural Science Foundation  of China.

\end{document}